\documentclass[11pt,reqno]{amsart}
\usepackage{amsmath,amstext,amsfonts,amsthm,amssymb}
\usepackage[all]{xy}
\usepackage{eucal}
\usepackage{graphicx}
\swapnumbers
\newtheorem{thm}[subsubsection]{Theorem}
\newtheorem{lem}[subsubsection]{Lemma}
\newtheorem{cor}[subsubsection]{Corollary}

{  \theoremstyle{definition}

}
\numberwithin{equation}{subsection}
\newcommand{\iso}{\overset{\sim}{\longrightarrow}}

\newcommand{\bea}{\begin{eqnarray*}}
\newcommand{\eea}{\end{eqnarray*}}
\newcommand{\bean}{\begin{eqnarray}}
\newcommand{\eean}{\end{eqnarray}}

\newcommand{\fa}{\mathfrak a}
\newcommand{\fg}{\mathfrak g}

\newcommand{\CA}{\mathcal{A}}

\newcommand{\CD}{\mathcal{D}}
\newcommand{\CF}{\mathcal{F}}

\newcommand{\CM}{\mathcal{M}}

\newcommand{\CO}{\mathcal{O}}

\newcommand{\CS}{\mathcal{S}}

\newcommand{\BH}{\mathbb{H}}
\newcommand{\BC}{\mathbb{C}}

\newcommand{\BZ}{\mathbb{Z}}

\newcommand{\nc}{\newcommand}

\nc{\Id}{\text{Id}}
\nc{\la}{\lambda}

\begin{document}


\newcommand{\Real}{\mathbb R}
\newcommand{\HH}{\mathbb H}
\newcommand{\QQ}{\mathbb Q}
\newcommand{\ZZ}{\mathbb Z}
\newcommand{\LL}{\mathbb L}
\newcommand{\VV}{\mathbb V}
\newcommand{\MM}{\mathbb M}
\newcommand{\PP}{\mathbb P}
\newcommand{\RR}{\mathbb R}

\newcommand{\cA}{\mathcal{A}}
\newcommand{\cB}{\mathcal{B}}
\newcommand{\cC}{\mathcal{C}}
\newcommand{\cext}{\cC ext}
\newcommand{\cD}{{\mathcal{D}}}
\newcommand{\cE}{{\mathcal{E}}}

\newcommand{\cF}{{\mathcal{F}}}
\newcommand{\cG}{{\mathcal{G}}}
\newcommand{\cH}{{\mathcal{H}}}
\newcommand{\cJ}{{\mathcal{J}}}

\newcommand{\cL}{{\mathcal{L}}}
\newcommand{\cM}{{\mathcal{M}}}
\newcommand{\cN}{\mathcal{N}}
\newcommand{\cO}{{\mathcal{O}}}
\newcommand{\cP}{{\mathcal{P}}}
\newcommand{\cQ}{{\mathcal{Q}}}
\newcommand{\cR}{{\mathcal{R}}}
\newcommand{\cS}{{\mathcal{S}}}
\newcommand{\cT}{{\mathcal{T}}}
\newcommand{\cV}{{\mathcal{V}}}
\newcommand{\cW}{{\mathcal{W}}}
\newcommand{\vext}{\cV ext}

\newcommand{\cZ}{{\mathcal{Z}}}

\title{Chiral De Rham complex over locally complete intersections}

\author{ Fyodor Malikov \and Vadim Schechtman}
\maketitle

\centerline{{\it To Borya Feigin on his 60th birthday}}

\begin{abstract}
Given a locally complete intersection $X\hookrightarrow Y$ we define a version of a derived chiral De Rham complex, thereby
``chiralizing'' a result by Illusie and Bhatt. A similar construction attaches to a graded ring a dg vertex algebra, which we prove to be
 Morita equivalent to a dg algebra of differential operators. For example, the dg vertex algebra associated to a fat point, which also
 arises in the Landau-Ginzburg model, is shown to be derived rational.
 \end{abstract}

\section{introduction}
\label{intro}
Grothendieck has attached an algebra of differential operators, $\CD_A$, to any commutative ring $A$.
However, if $\text{Spec}A$ is not a smooth variety then this algebra may be pathological,
and the usual approach is to abandon the algebra and define the category of $\CD$-modules
by embedding $\text{Spec}A$ into a smooth algebraic variety $Y$ and letting $\CD-mod$ be the category of $\CD_Y$-modules supported on $\text{Spec}A$. (All algebras will be supposed to be 
 of finite type over $\BC$.) 

The definition of an algebra of {\em chiral} differential operators \cite{MSV, GMS} does not
even hint at a method to construct it on a singular variety. In this paper we suggest a 
homotopical version of such construction which uses the ideas 
of derived algebraic geometry.

Let $P$ be a {\em finite} polynomial differential graded (dg) algebra and $P\rightarrow A$ a dga resolution.
(All our dg algebras are nonpositively graded with differential of degree +1; $A$ is placed in degree 0 and has zero differential. 
As is well known, the existence of a finite dga resolution is equivalent to $A$ being a complete intersectiion.)  We  verify that $\CD_P$ is a dg algebra determined uniquely
up to homotopy, and so we denote by $L\CD_A$ its homotopy class. D.Gaitsgory has explained to
us that the derived category of dg $\CD_P$-modules is equivalent to the derived category of
$\CD$-modules over $\text{Spec}A$ defined as above.

It is perfectly clear what the dg algebra of chiral differential operators, $\CD^{ch}_P$, is, and our proposal is to consider the assignment $A\mapsto \CD^{ch}_P$, or its De Rham version
$A\mapsto \CD^{ch}_{\Omega^\bullet_P}$. The merit of this proposal depends on how natural
such 
assignments are, and we admit that we do not know much about this.  Let us now describe the results.

The fundamental problem is that there is no such thing as a free vertex algebra, and we have to rely
on stopgap measures. 
If we fix a locally complete intersection  $X\hookrightarrow Y$, where $Y$ is a smooth
algebraic variety, then the De Rham version of our proposal gives $L\Omega^{ch}_{X\hookrightarrow Y}$, a sheaf  on $X$ with values in an appropriate homotopy category
of vertex algebras.  

The construction manifestly depends on the embedding, but part of it does not. The sheaf
$L\Omega^{ch}_{X\hookrightarrow Y}$ is conformally graded, and its 0th component, $L\Omega^{ch}_{X\hookrightarrow Y}[0]$,
is essentially the derived De Rham complex, $L\Omega^\bullet_X$, of Illusie \cite{Ill2}.   There arises a natural projection 
$L\Omega^{ch}_{X\hookrightarrow Y}\rightarrow L\Omega^{ch}_{X\hookrightarrow Y}[0]$, which  gives an isomorphism
$\BH^\bullet(X,L\Omega^{ch}_{X\hookrightarrow Y})\iso \BH^\bullet(X,L\Omega^{ch}_{X\hookrightarrow Y}[0])$.

On the other hand, Grothendieck's classic  result on the algebraic description of singular cohomology has been generalized by Illusie in the case
of locally complete intersections and Bhatt \cite{Bh} in general as follows: $\BH^\bullet(X,L\Omega^\bullet_X)\iso H^\bullet(X,\BC)$. Combined the last two isomorphisms give
\[
\BH^\bullet(X,L\Omega^{ch}_X)\iso H^\bullet(X,\BC).
\]
This can be thought of as a {\em chiralization} of the Bhatt-Illusie theorem. Of course, this isomorphism is also a generalization of an analogous
result for the chiral De Rham complex $\Omega^{ch}_X$ over a smooth $X$ proved in \cite{MSV}.

The homotopy content of the assignment $A\mapsto \CD^{ch}_P$ is mirkier, although we do prove that given two polynomial dga resolutions
$P\rightarrow A\leftarrow Q$ the dg vertex algebras $\CD^{ch}_P$ and $\CD^{ch}_Q$ are quasiisomorphic  provided $A,P,Q$ all carry an ``inner''
grading. More interesting insight into the situation is provided by a representation-theoretic point of view.  We show that the categories of dg modules
over $\CD_P$ and $\CD_P^{ch}$ are equivalent, which via the 
Gaitsgory's result mentioned above implies that the derived category of $\CD$-modules
over $\text{Spec}A$ is equivalent to the derived category of dg modules over $\CD^{ch}_P$. This result is of course independent of the choice of $P$.

A nice illustration is provided by $A=\BC[x]/(x^n)$.  In this case $P$ is naturally chosen to be the Koszul resolution $K(\BC[x],x^n)$, and the dg algebra
$\CD_{K(\BC[x],x^n)}$ is quasiisomorphic to the matrix algebra $gl_n(\BC)$. It follows that in this case the derived category of dg modules over 
$\CD^{ch}_{K(\BC[x],x^n)}$
is equivalent to the derived category of finite dimensional vector spaces. One may want to think of this result as either a version of Morita equivalence
or an ``odd Stone-von Neumann theorem.''

The Milnor algebra $\BC[x]/(x^n)$ was for us a starting point.  The dg vertex algebra $\CD^{ch}_{K(\BC[x],x^n)}$ is intimately  related to the 
Landau-Ginzburg model \cite{Witt}, and its cohomology was computed in a beautiful paper by B.Feigin and A.Semikhatov \cite{FS}. A little thought
will convince the reader that the above description of the derived category of dg modules can be given the following form: {\em the Feigin-Semikhatov
vertex algebra is derived rational.} In fact, whenever $\text{Spec}A$  is a union of finitely many fat points, the corresponding vertex  algebra is
derived rational. 

Another motivational example was provided by the Berkovits model, see e.g. \cite{BN}.  It was explained by N.Nekrasov \cite{N} that the Berkovits
model requires an algebra of chiral differential operators on the space of pure spinors, which is a cone singular at its vertex. This space, however,
is not a complete intersection, and dealing with infinite resolutions appears to be a major challenge in the present context.

{\em Acknowledgments.}  The authors have greatly benefited from conversations with  A.Beilinson, D.Gaitsgory, A.Gerasimov, V.Gorbounov, V.Hinich,
M.Kapranov, M.Kontsevich, N.Nekrasov.  Parts of this work were done when the authors were visiting MPIM in Bonn and (the first author) IHES. We are grateful to these institutions for
excellent working conditions. F.M. was partially supported by an NSF grant.

\section{derived chiral de rham complex}
\label{derived chiral de rham complex}
\subsection{Vertex algebras.}  We shall merely set up the notation referring the reader to sources such as \cite{FBZ,K,BD} for the introduction to the subject.

\subsubsection{ } 
\label{list-vert-alg-notat}
A vertex algebra is a vector  super-space $V$ carrying a countable family
of multiplication
\[
_{(n)}: V\times V\rightarrow V,\; a\otimes b\mapsto a_{(n)}b, n\in\BZ,
\]
and a derivation
\[
T: V\rightarrow V\text { s.t. }\; T(a_{(n)}b)=(T a)_{(n)}b+a_{(n)}(T b).
\]
One way to state the axioms that these multiplications satisfy is to form {\em fields}
\[
V\ni a\mapsto a(z)=\sum_{n\in\BZ}a_{(n)}z^{-n-1}
\]
and explain in what sense
\[
a(z)b(w)=b(w)a(z)=\sum_{n\gg -\infty}\frac{(a_{(n)}b)(w)}{(z-w)^{-n-1}}.
\]
A conformal grading on $V$ is a direct sum decomposition 
\[
V=\bigoplus_{n=0}^{\infty} V[n]\text{ s.t. } V[n]_{(i)}V[m]\subset V[n+m-i-1],
\]
and $a\in V[n]$ is said to have conformal weight $n$.

An element of $\text{End}_{\BC}(V)$ is called a {\em derivation} if it is a derivation of each of the products. For example, for each $a\in V$, $a_{(0)}$
(regarded as a map $v\mapsto a_{(0)}v$ is a derivation. The space of all derivations, $\text{Der}(V)$, is a Lie (super)algebra.

An element $\delta\in\text{Der}(V)$ is called a differential if its square is 0.  If, in additiion, $V$ is a graded vector space (typically this grading is different from the conformal
one) s.t. $\delta$ has degree 1, then the pair $(V,\delta)$ is called a differential graded vertex algebra.

\subsubsection{ } 
\label{local-defn-of-cdo}Let $X$ be a smooth affine variety over $\BC$ and denote by $A=\BC[X]$  its coordinate ring.  A  coordinate system, $\vec{x}$,  is
is a collection of  elements $x_1,...,x_N\in A$,
$\partial_{x_1},...,\partial_{x_N}\in\text{Der}(A)$ s.t.  $\partial_{x_i}(x_j)=\delta_{ij}$ and the set $\{\partial_1,\partial_2,...\}$ is a basis (over $A$) of $\text{Der}(A)$.  One has:
$\{dx_1,...,dx_N\}$ is a basis of $\Omega_A$ and $[\partial_i,\partial_j]=0$.

If exists, a  coordinate system defines a  vertex algebra structure on $\BC[J_\infty TX]$, the coordinate ring of the space of jets in the cotangent bundle $T^*X$.
Namely,   $\vec{x}$ gives an identification
$\BC[J_\infty T^*X]\stackrel{\sim}{\rightarrow}A[x_i^{(n)},\partial_{x_j}^{(m)}; \;1\leq i,j\leq \text{dim} A, m,n\geq 0]$. The latter space is an algebra with derivation
$\nabla$ s.t. $\nabla(a^{(n)})=a^{(n+1)}$ for $a$ either $x_i$ or $\partial_{x_j}$, and it carries a vertex algebra structure determined by the conditions
\begin{eqnarray}
(\partial_{x_i})_{(0)}f&=&\partial_{x_i}(f),\; (\partial_{x_i})_{(m)}f=0\text{ if }m>0, \forall f\in A,
\label{defn-cdo-1}\\
(x_i)_{(-n-1)} P&=&x_i^{(n)}\cdot P, (\partial_{x_j})_{(-n-1)}P=\partial_{x_j}^{(n)}\cdot P\text{ if }n\geq0,\; \forall P\in\BC[J_\infty TX],
\label{defn-cdo-2}\\
T&=&\nabla.
\label{defn-cdo-3}
\end{eqnarray}
Here and elsewhere $a$ is identified with $a^{(0)}$ thus encoding the canonical embeddings $\BC[X]\hookrightarrow\BC[T^*X]\hookrightarrow\BC[J_\infty T*X]$.

This vertex algebra is denoted $\CD^{ch}_{A,\vec{x}}$ (or $\CD^{ch}_{X,\vec{x}}$ depending on the context) and called an {\em algebra of chiral differential operators }over $A$ (or $X$) for the reason
that (\ref{defn-cdo-1}) resemble (and extend) the familiar $[\partial_{x_i},f]=\partial_{x_i}(f)$.  

This construction has  superalgebra extensions ($\CD^{ch}_A$ is of course assumed purely even), the most important of which occurs when $A$ is
replaced with the De Rham algebra $\Omega^\bullet_A$; or, equivalently,  $X$ is replaced with
the (total space of the) supervector bundle  $\Pi TX$, where $\Pi TX$ is $TX$ with fibers of the projection $TX\rightarrow X$ declared odd.  In the presence of  coordinates
the passage from $X$ to $\Pi TX$ consists in adjoining Grassman variables $\phi_i$, $1\leq i\leq \text{dim} A$, and the jet-algebra becomes
\[
\BC[J_\infty T^*(\Pi TX)]\stackrel{\sim}{\rightarrow}\BC[J_\infty TX]\otimes\BC[\phi_i^{(n)},\partial_{\phi_i}^{(n)};\; 1\leq i\leq \text{dim} A, n\geq 0].
\]
\begin{sloppypar}
The vertex algebra structure on $\BC[J_\infty T^*(\Pi TX)]$ is defined by an obvious extension of (\ref{defn-cdo-1},\ref{defn-cdo-2},\ref{defn-cdo-3}), where we
allow  (1) $f$ to be an element of $\BC[\Pi TX]$, and (2)  to replace $x_i$ ($\partial_{x_i}$ resp.) with $\phi_i$ ($\partial_{\phi_i}$ resp.) The result is denoted
$\CD^{ch}_{\Omega^\bullet_A,\vec{x}}$and is an example of a superalgebra of chiral differential operators. 

Both $\CD^{ch}_{A,\vec{x}}$ and $\CD^{ch}_{\Omega^\bullet_A,\vec{x}}$ are 
  conformally graded s.t. $D^{ch}_{A,\vec{x}}[0]=A$ and $D^{ch}_{\Omega_A^\bullet,\vec{x}}[0]=\Omega_A^\bullet$, and this condition uniquely determines the grading.
  \end{sloppypar}
  
  \subsubsection{ } The inquisitive  reader must have  noticed that not only $\CD^{ch}_{\Omega^\bullet_A,\vec{x}}$ depends on a choice of a coordinate system, it is not related to the tangent bundle (in $\Pi TX$) in
  any serious way, since a trivialization of any  trivial vector bundle $E\rightarrow X$ will produce in the same way the same vertex algebra. However,  further pursuing the
  De Rham interpretation is fruitful. 
  
  To begin with, just as $\Omega^\bullet_A$ is a differential graded algebra with differential $d_{DR}=\sum_i\phi_i\partial_{x_i}$,
  $\CD^{ch}_{\Omega^\bullet_A}$ is a differential graded vertex algebra with differential $d_{DR}^{ch}=\sum_i((\phi_i)_{(-1)}\partial_{x_i})_{(0)}$.
  (The grading on the former is defined by assigning $\phi\mapsto1$, one on the latter by assigning $\phi\mapsto1, \partial_{\phi}\mapsto-1$, which gives $d_{DR}$
  and $d_{DR}^{ch}$ degree 1.)

 The crux of the matter is that the assignment $A\mapsto(\CD^{ch}_{\Omega^\bullet_A,\vec{x}},d^{ch}_{DR})$ is functorial w.r.t. isomorphisms; let us explain this.  Any isomorphism $f: A\rightarrow B$ defines an isomorphism $df:\Omega^\bullet_A\rightarrow\Omega^\bullet_B$, hence  an isomorphism of the
  corresponding jet-spaces  $J_\infty T^*(\Pi TX_B)\rightarrow J_\infty T^*(\Pi TX_A)$, hence  an isomorphism $df_\infty: \BC[J_\infty T^*(\Pi TX_A)]\rightarrow
   \BC[J_\infty T^*(\Pi TX_B)]$. In the presence of coordinates, this tautological isomorphism is defined by some explicit and  classical formulas; e.g.,   the induced map of vector
   fields $df_\infty: \text{Der}(\Omega_A)\rightarrow\text{Der}(\Omega_B)$ is defined by $df_\infty(\xi)(b)=df(\xi df^{-1}(b))$ and in terms of $\vec{x}\subset A$, 
   $\vec{\tilde{x}}\subset B$ becomes
   \begin{equation}
   \label{map-on-jets-DR}
   df_\infty(\partial_{x_i})=f(\partial_{x_i}f^{-1}(\tilde{x}_\alpha))\partial_{\tilde{x}_\alpha}+
   f(\partial_{x_i}\partial_{x_j}(f^{-1}\tilde{x}_\alpha))\partial_{\tilde{x}_\beta}(fx_j)\tilde{\phi_\beta}\partial_{\tilde{\phi}_\alpha}.
   \end{equation}
   Although $\CD^{ch}_{\Omega^\bullet_A,\vec{x}}$ is identified with $\BC[J_\infty T^*(\Pi TX_A)]$, the induced map $df_\infty:\CD^{ch}_{\Omega^\bullet_{A,\vec{x}}}\rightarrow
  \CD^{ch}_{\Omega^\bullet_{B,\vec{\tilde{x}}}}$ is not a vertex algebra morphism. It can, however, be corrected to become one.
  One has \cite{MSV}
  
  \begin{lem}
  \label{lift-to-chiral-dr}
  Let $f:A\rightarrow B$ be an algebra isomorphism. There is a unique differential graded
   vertex algebra isomorphism $df^{ch}_{\vec{x},\vec{\tilde{x}}}:\CD^{ch}_{\Omega^\bullet_{A,\vec{x}}}\rightarrow
  \CD^{ch}_{\Omega^\bullet_{B,\vec{\tilde{x}}}}$ s.t. (cf. (\ref{map-on-jets-DR}) )
  \begin{eqnarray}
  A\ni a&\mapsto& f(a)\nonumber\\
  \phi_i&\mapsto& \partial_{\tilde{x}_\alpha}(fx_i)\phi_\alpha\nonumber\\
  \partial_{\phi_i}&\mapsto&f(\partial_{x_\alpha}(f^{-1}\tilde{x}_\alpha))\partial_{\tilde{\phi}_\alpha}\nonumber\\
\partial_{x_i}&\mapsto&f(\partial_{x_i}f^{-1}(\tilde{x}_\alpha))\partial_{\tilde{x}_\alpha}+
   f(\partial_{x_i}\partial_{x_j}(f^{-1}\tilde{x}_\alpha))\partial_{\tilde{x}_\beta}(fx_j)\tilde{\phi_\beta}\partial_{\tilde{\phi}_\alpha}.\nonumber
  \end{eqnarray}
  These isomorphisms satisfy the associativity condition: for any isomorphisms $f:A\rightarrow B$, $g:B\rightarrow C$ and bases $\vec{x}\subset A$, $\vec{y}\subset B$,
  $\vec{z}\subset C$
  \[
  d(g\circ f)^{ch}_{\vec{x},\vec{z}}=dg^{ch}_{\vec{y},\vec{z}}\circ df^{ch}_{\vec{x},\vec{y}}.
  \]
  Furthermore, each $df^{ch}_{\vec{x},\vec{\tilde{x}}}$ preserves the conformal grading mentioned at the end of sect.~\ref{local-defn-of-cdo}.
  \end{lem}
  
  Attached to $A$ is a collection of vertex algebras $\{\CD^{ch}_{\Omega^\bullet_A,\vec{x}},\;\vec{x}\subset A\}$, each pair among which is connected by an isomorphism
  $d(Id)^{ch}_{\vec{x},\vec{y}}:\CD^{ch}_{\Omega^\bullet_A,\vec{x}}\rightarrow\CD^{ch}_{\Omega^\bullet_A,\vec{y}}$. These isomorphisms satisfy 
  $d(Id)^{ch}_{\vec{x},\vec{z}}=d(Id)^{ch}_{\vec{y},\vec{z}}\circ d(Id)^{ch}_{\vec{x},\vec{y}}$. This means that there is a well-defined vertex algebra, $\CD^{ch}_{\Omega^\bullet_A}$,
  attached to $A$. Furthermore, the same associativity property implies that a collection of isomorphisms 
  $df^{ch}_{\vec{x},\vec{y}}:\CD^{ch}_{\Omega^\bullet_A,\vec{x}}\rightarrow\CD^{ch}_{\Omega^\bullet_B,\vec{y}}$ gives a well-defined isomorphism
$df^{ch}:\CD^{ch}_{\Omega^\bullet_A}\rightarrow\CD^{ch}_{\Omega^\bullet_B}$. This proves

\begin{cor}
\label{ch-dr-functor-cor}
The assignments $A\mapsto \CD^{ch}_{\Omega^\bullet_A}$, $f\mapsto df^{ch}$ define a functor from the category of $\BC$-algebras that admit
a coordinate system and isomorphisms to the category of dg vertex algebras.
\end{cor}
  
  In particular,  there is a natural group homomorphism 
  \[
  \text{Aut}(A)\rightarrow \text{Aut}(\CD^{ch}_{\Omega^\bullet_A}), \; f\mapsto df^{ch}.
  \]
  This result has a simpler  infinitesimal analogue (recorded e.g. in \cite{MS}):
  \begin{lem}
  \label{functor-rings-dgva-rev-infinit-incarn}
  The assignment
  \[
  f_i\partial_{x_i}\mapsto ((\partial_{x_i})_{(-1)}f_i-(\partial_{\phi_i})_{(-1)}\phi_\alpha\partial_{x_\alpha}(f_i))_{(0)}
  \]
  defines a natural dg Lie algebra morphism
  \[
  \text{Der}(A)\rightarrow\text{Der}(\CD^{ch}_{\Omega^\bullet_A}).
  \]
  \end{lem}

\subsubsection{ }  
\label{defn-of-disting-chiral-de-rham}
Let $A$ be an algebra that admits a coordinate system and $S\subset A$ be multiplicative. Since all the maps involved have been defined using
various differential operators, one obtains a  vertex dg algebra structure on the localization $(\CD^{ch}_{\Omega^\bullet_A})_S$, which we prefer to denote by
$\CD^{ch}_{\Omega^\bullet_{A_S}}$, along with a vertex dg algebra morphism $\CD^{ch}_{\Omega^\bullet_A}\rightarrow \CD^{ch}_{\Omega^\bullet_{A_S}}$.
Corollary~\ref{ch-dr-functor-cor} now implies, by glueing, that  over each smooth algebraic variety $X$ there is a sheaf of vertex differential graded
algebras $\CD^{ch}_X$, and this sheaf is natural in that for any isomorphism $f: X\rightarrow Y$  there is a natural isomorphism 
$\CD^{ch}_X\rightarrow f^{-1}\CD^{ch}_Y$. This particular algebra of chiral differential operators will be denoted by $\Omega^{ch}_X$ or $\Omega^{ch}_A$ if
$X=\text{Spec}(A)$.

\subsection {A derived version}

\subsubsection{ }
\label{defn-of-alg-homotop}
Let $R$ be a $\BC$-algebra. A dg $R$-algebra will mean one with grading bounded from the right and differential $\partial$ of degree +1.  

Given a dg algebra morphism $f: A\rightarrow B$, define $\text{Der}_f(A,B)$ to be the space of $R$-linear maps $\phi:A\rightarrow B$ that satisfy the Leibniz condition:
$\phi(xy)=\phi(x)f(y)+(-1)^{\tilde{\phi}\tilde{x}}f(x)\phi(y)$.

Two dg algebra homomorphisms $f,g: A\rightarrow B$ are homotopic if there is a polynomial family of dg algebra morphisms
\[
h_t: A\rightarrow B\text{ s.t. } h_{t=0}=f, h_{t=1}=g,
\]
and a polynomial family of degree -1 $R$-morphisms
\[
\alpha_t: A\rightarrow B\text{ s.t. }\alpha_{t=t_0}\in\text{Der}_{h_{t=t_0}}(A,B), [\partial,\alpha_{t=t_0}]=\frac{d}{dt}|_{t=t_0}h_t\;\forall t_0\in\BC.
\]
This concept is what is known in the closed model category 
theory as the right homotopy with  $B\otimes_R\Omega_{\BC}$ chosen as a path object. As defined, homotopy is
not necessarily an equivalence relation, but it is if the algebras $A$ and $B$ are semi-free
\footnote{Recall that a dga algebra is called semi-free if it is a polynomial ring and its set of variables carries an exhaustive  filtration  $X_0=\emptyset\subset X_1\subset X_2\subset\cdots$ s.t.  $\partial(I_{n+1})\subset\BC[I_n]$.} (or, more generally, cofibrant.)

\subsubsection{ }  
\label{derived-version-of-a-ring}
Let  $R$ be a coordinate ring of a smooth affine variety over $\BC$, $J\subset R$ a regular ideal, and $A=R/J$. A choice of a regular generating
system $\{f_1,...,f_r\}\subset J$ defines a Koszul resolution, i.e., a dg algebra $(R[\xi_1,...,\xi_r],\partial)$, where $R$ is declared even of degree 0,  $\xi$'s  odd
of degree -1,  $\partial=f_i\partial_{\xi_i}$
and a quasiisomorphism $R[\xi_1,...,\xi_r]\rightarrow A$ defined to be the quotienting out by the ideal generated by $J$ and $\{\xi_1,...,\xi_r\}$.  Denote this dg algebra by
$K(R,\vec{f})$.

Another choice of generators, $\{g_1,...,g_r\}$, defines another Koszul resolution, $K(R,\vec{g})$. There are $r\times r$ matrices $F$ and $G$ with coefficients in $R$ s.t.
$f_i=G_{i\alpha}g_\alpha$, $g_i= F_{i\alpha}f_\alpha$.  This gives 2 dg $R$-algebra isomorphisms
\begin{equation}
\label{iso-koszul-dga}
\tilde{G}:K(R,\vec{f})\rightarrow K(R,\vec{g}),\; \xi_i\mapsto G_{i\alpha}\eta_\alpha;\;
\tilde{F}: K(R,\vec{g})\rightarrow K(R,\vec{f}),\; \eta_i\mapsto F_{i\alpha}\xi_\alpha.
\end{equation}
The matrices $F$ and $G$ are not unique (but they are modulo $J$, since $J/J^2$ is a free $A$-module), and so the compositions $\tilde{F}\circ\tilde{G}$, 
$\tilde{G}\circ\tilde{F}$ do not have to be an identity, but they are up to homotopy.  This follows from the closed model category theory on the grounds that 
both $K(R,\vec{f})\rightarrow A$ and $K(R,\vec{g})\rightarrow A$ are cofibrant resolutions in the category of dg $R$-algebras. Therefore, the assignment $A\mapsto K(R,\vec{f})$
is a functor from the category of $R$-algebras and isomorphisms to the homotopy category of dg $R$-algebras.

\subsubsection{ }
\label{derived-version-of-chiral-ring}
We would like to ``chiralize'' the considerations of sect.~\ref{derived-version-of-a-ring} in order to attach a dg vertex algebra to the algebra $A$.  Consider a Koszul
resolution $K(R,\vec{f})\rightarrow A$.  All the considerations that led in sect.~\ref{defn-of-disting-chiral-de-rham} to the chiral De Rham complex $\Omega^{ch}_R$
go through unchanged in the case of a regular supercommutative $\BC$-algebra, such as $R[\xi_1,...,\xi_r]$. This amounts to adjoining to $\Omega^{ch}_R$ 4 groups of
variables: $\xi_i, \partial_{\xi_i}$, all odd, and $\xi_i^*,\partial_{\xi^*_i}$, all even, $\xi_i^*$ having the meaning of $d_{DR}\xi_i$,
and stipulating the obvious versions of the relations from sect.~\ref{local-defn-of-cdo}.  For the reasons that will become apparent later, it is convenient to complete this vertex
algebra by allowing power series in $\xi_i^*$, $1\leq i\leq r$: it was one of observations made in \cite{MSV} that the vertex algebra multiplications extend to such completions by continuity.

The result is a not simply a dg but a bi-differential bi-graded
vertex algebra $\Omega^{ch}_{K(R,\vec{f})}$,  since in addition to $d^{ch}_{DR}$ we now have $\partial^{ch}$ defined to be 
\[
\partial^{ch}=(f_i\partial_{\xi_i}-(d_{DR}f_i)\partial_{\xi_i^*})_{(0)},\text{ where } d_{DR}f=\phi_i\partial_{x_i}(f),
\]
 which is a natural extension
of the Koszul differential $\partial=f\partial_{\xi_i}$, and one readily verifies that $[d^{ch}_{DR},\partial^{ch}]=0$. The corresponding bi-grading where the degrees
of $\partial^{ch}$ and $d^{ch}_{DR}$ are $(-1,0)$ and $(0,1)$ resp. is defined by assigning degree $(-1,0)$ to $\xi_i$'s, $(0,1)$ to $\phi_i$'s, $(-1,1)$ to $\xi_i^*$'s,
$(0,0)$ to elements of $R$, and finally, minus that to the corresponding vector fields $\partial_{\xi_i}$, etc.

Lemma~\ref{lift-to-chiral-dr} carries over to the present situation and implies that isomorphisms (\ref{iso-koszul-dga}) lift to isomorphisms if bi-differential bi-graded vertex
algebras
\begin{equation}
\label{dg-chiral-koszul-iso}
G^{ch}:\Omega^{ch}_{K(R,\vec{f})}\stackrel{\longrightarrow}{\longleftarrow}\Omega^{ch}_{K(R,\vec{g})}: F^{ch}.
\end{equation}
Neither $F^{ch}\circ G^{ch}$ nor $G^{ch}\circ F^{ch}$ have to be the identity, but they are homotopic to identity. Let us explain this.

\subsubsection{ }
\label{homo-for-vert-alg}
The discussion in sect.~\ref{defn-of-alg-homotop} is easily adjusted to the vertex algebra case.  Given a dg vertex 
algebra morphism $f: A\rightarrow B$, define $\text{Der}_f(A,B)$ to be the space of $\BC$-linear maps $\phi:A\rightarrow B$ that satisfy the Leibniz condition for each of the 
multiplications:
$\phi(x_{(n)}y)=\phi(x)_{(n)}f(y)+(-1)^{\tilde{\phi}\tilde{x}}f(x)_{(n)}\phi(y)$, $n\in \BZ$.

Two dg vertex algebra homomorphisms $f,g: A\rightarrow B$ are homotopic if there is a polynomial family of dg vertex algebra morphisms
\[
h_t: A\rightarrow B\text{ s.t. } h_{t=0}=f, h_{t=1}=g,
\]
and a polynomial family of degree -1 morphisms
\[
\alpha_t: A\rightarrow B\text{ s.t. }\alpha_{t=t_0}\in\text{Der}_{h_{t=t_0}}(A,B), [\partial,\alpha_{t=t_0}]=\frac{d}{dt}|_{t=t_0}h_t\;\forall t_0\in\BC.
\]
Creating a version of closed model category theory for vertex 
algebras is problematic because of  the lack of free vertex algebras, and the proposed definition may be of
limited use, but in the specific situation of sect.~\ref{derived-version-of-chiral-ring} there is a canonical lift of a homotopy for polynomial rings to that for the corresponding
vertex algebra as follows.  

The set-up  is apparently  this: a dg algebra $K(R,\vec{f})$, a dg vertex algebra $\Omega^{ch}_{K(R,\vec{f})}$ with differential $\partial^{ch}$,
 an isomorphism $g: K(R,\vec{f})\rightarrow K(R,\vec{f})$ and a homotopy $(h_t,\alpha_t)$ of $g$ and $Id$. We have a lift of
$g$ to $g^{ch}\in\text{Aut}(\Omega^{ch}_{K(R,\vec{f})})$, and we want to find a homotopy of $g^{ch}$ and $Id$. 

 Since, by definition,
$\alpha_{t=t_0}\in\text{Der}_{h_{t=t_0}}(K(R,\vec{f}))$,  the composition $h^{-1}_{t=t_0}\circ\alpha_{t=t_0}\in\text{Der}(K(R,\vec{f}))$.  Lemma~\ref{functor-rings-dgva-rev-infinit-incarn} has an obvious extension in the case where $A$ is a supercommutative ring (some signs must be changed) and implies that  there are lifts 
of $h_{t=t_0}$ to $h^{ch}_{t=t_0}\in\text{Aut}(\Omega^{ch}_{K(R,\vec{f})})$ and  of $h^{-1}_{t=t_0}\circ\alpha_{t=t_0}\in\text{Der}(K(R,\vec{f}))$ to $(h^{-1}_{t=t_0}\circ\alpha_{t=t_0})^{ch}\in\text{Der}(\Omega^{ch}_{K(R,\vec{f})})$.  Lemma~\ref{functor-rings-dgva-rev-infinit-incarn} applied one more time gives: the pair $(h_t^{ch},\alpha^{ch}_t)$
with $\alpha^{ch}_{t=t_0}$ defined to be $h_{t=t_0}^{ch}\circ(h^{-1}_{t=t_0}\circ\alpha_{t=t_0})^{ch}$ defines a homotopy of $g^{ch},\text{Id}\in\text{Aut}(\Omega^{ch}_{K(R,\vec{f})})$.

Notice that (1) by construction these homotopies of ``classical origin'' indeed define an equivalence relation, and (2) the constructed isomorphisms preserve the 2 differentials
($\partial^{ch}$ and $d_{DR}^{ch}$) and the 2 gradings.

To summarize:

\begin{lem}
\label{derived-chir-de-rham-ring}
Under the assumptions  of sect.~\ref{derived-version-of-a-ring},  the bi-differential bi-graded vertex algebra $\Omega^{ch}_{K(R,\vec{f})}$ is determined uniquely
up to homotopy.
\end{lem}
We shall denote this vertex algebra by $L\Omega^{ch}_{R\mapsto A}$ and call it the derived chiral De Rham complex over $A$ relative to the projection $R\rightarrow A$.

\subsubsection{ }
\label{extens-to-varieties}
Lemma~\ref{derived-chir-de-rham-ring} implies, by gluing, that if $X$ is a smooth algebraic variety and $Y\hookrightarrow X$ is a locally complete intersection, then $Y$ carries 
a sheaf with values in the homotopy category of bi-differential bi-graded vertex algebras  to be denoted $L\Omega^{ch}_{Y\hookrightarrow X}$ and called
the derived chiral De Rham complex of $Y$ relative to $Y\hookrightarrow X$.  Indeed, each point in $Y$ has an affine neighborhood $U$ that is a complete intersection. For each
affine $V\subset X$ s.t. $V\cap Y=U$ consider $L\Omega^{ch}_{\CO_X(V)\rightarrow\CO_Y(U)}$.  
It follows from sect.~\ref{defn-of-disting-chiral-de-rham} that the family $\{L\Omega^{ch}_{\CO_X(V)\rightarrow\CO_Y(U)};\; V \text{ s.t. }V\cap Y=U\}$ is an inductive system.
The sheaf  $L\Omega^{ch}_{Y\hookrightarrow X}$ is a unique sheaf s.t. for
each such  $U$
\[
L\Omega^{ch}_{Y\hookrightarrow X}(U)=\lim_{ V: \stackrel{\longrightarrow}{V\cap Y}=U}L\Omega^{ch}_{\CO_X(V)\rightarrow\CO_Y(U)}.
\]
In addition to being bi-graded (as in ``bi-graded bi-differential''), $L\Omega^{ch}_{Y\hookrightarrow X}$ is conformally graded:
\[
L\Omega^{ch}_{Y\hookrightarrow X}=\bigoplus_{n=0}^{+\infty}L\Omega^{ch}_{Y\hookrightarrow X}[n],
\]
as all the manipulations have preserved the
conformal grading, cf. Lemma~\ref{lift-to-chiral-dr}. Notice that each differential has conformal weight 0, and so each component $L\Omega^{ch}_{Y\hookrightarrow X}[n]$ is a bicomplex. 

Now consider $L\Omega^{ch}_{Y\hookrightarrow X}$ as a total complex.  There is an obvious projection
\[
L\Omega^{ch}_{Y\hookrightarrow X}\rightarrow L\Omega^{ch}_{Y\hookrightarrow X}[0],
\]
and one observes that $L\Omega^{ch}_{Y\hookrightarrow X}[0]$ is essentially $L\Omega^{\bullet}_Y$, the derived De Rham complex of Illusie \cite{Ill2}\footnote{ It is true that Illusie works
with simplicial rather than dg commutative algebras, but these categories are related by a Quillen
equivalence $(N,N^*)$, see \cite{Q}, Remark at the bottom of p.223 (also \cite{SchSh} p.289). In characteristic 0 this gives an equivalence of both the approaches. We leave the details out.}.
At least there is a morphism
\[
L\Omega^{\bullet}_Y\rightarrow L\Omega^{ch}_{Y\hookrightarrow X}[0],
\]
which is easily seen to be a quasiisomorphism. A passage to hypercohomology gives a morphism
\[
\BH(Y,L\Omega^{ch}_{Y\hookrightarrow X})\rightarrow\BH(Y,L\Omega^\bullet_Y).
\]
On the other hand, Illusie \cite{Ill2} in the case of a locally complete intersection and Bhatt \cite{Bh} in general has proved an isomorphism
\[
\BH(Y,L\Omega^\bullet_Y)\iso H(Y,\BC),
\]
thus generalzing Grothendieck's algebraic description of singular cohomology.
\begin{thm}
\label{chiral-bhatt-thm}
The composite map
\[
\BH(Y,L\Omega^{ch}_{Y\hookrightarrow X})\rightarrow\BH(Y,L\Omega^\bullet_Y)\iso H(Y,\BC)
\]
is an isomorphism.
\end{thm}

{\em Proof.} It suffices to show that $L\Omega^{ch}_{Y\hookrightarrow X}\rightarrow L\Omega^{ch}_{Y\hookrightarrow X}[0]$ is a quasiisomorphism.
The question is local, and we will use local coordinates to show that the restriction of $d=d^{ch}_{DR}+\partial^{ch}$ to each nonzero conformal weight component is
homotopic to 0. We introduce, as in \cite{MSV}, the $N=2$ superconformal conformal algebra generators 
\[
L= Tx_i \partial_{x_i}+T\phi_i\partial_{\phi_i}+T\xi_i\partial_{\xi_i}+T\xi^*_i\partial_{\xi^*_i}
\]
and 
\[
G= Tx_i\partial_{\phi_i}+T\xi_i\partial_{\xi^*_i}
\]

 and compute (using the Wick theorem) to the effect that
\[
[G_{(1)},d^{ch}_{DR}]=L_{(1)},
\]
\[
[G_{(1)},\partial^{ch}]=0,
\]
and so
\[
[G_{(1)},d]=L_{(1)}.
\]
The first of these is familiar and was used in \cite{MSV} (without $\xi_i$ and $\xi_i^*$) to prove a similar statement in the smooth case
(i.e. without  $\xi_i$,  $\xi_i^*$, and $\partial^{ch}$).  As to the 2nd, one has
\[
[G_{(1)},\partial^{ch}]=((f_i\partial_{\xi_i}-(d_{DR}f_i)\partial_{\xi_i^*})_{(0)}(Tx_i\partial_{\phi_i}+T\xi_i\partial_{\xi^*_i}))_{(1)}=(Tf_i\partial_{\xi_i^*}-Tf_i\partial_{\xi_i^*})_{(1)}=0.
\]

Finally, observing that $L_{(1)}$ is the operator
that defines conformal grading, we deduce from the last equation that  $G_{(1)}/n$ is a homotopy of $d$ and 0 when restricted to $L\Omega^{ch}_{Y\hookrightarrow X}[n]$,
$n\not=0$.

\section{derived algebras of chiral and ordinary differential operators:  morita equivalence}
\label{derived algebras of chiral and ordinary differential operators; a morita duality}

\subsection{Derived algebras of differential operators.}
\label{Derived algebras of differential operators.} Let $A$ be a finitely generated commutative
 $\BC$-algebra, which we will regard as a differential graded algebra with  differential 0 and placed in
(cohomological) degree 0.  Assume further that there exists a dga resolution $P\rightarrow A$ s.t. $P$  is a {\em finitely generated} polynomial dga with generators placed in nonpositive degrees and
differential $\partial\subset\text{Der}(P)$. Now form the dg algebra of differential operators $\CD_{P}$ with differential $[\partial,.]$ and grading defined by $\text{deg}\partial_x=-\text{deg}x$ for each
generator $x$. 

\begin{lem}
\label{coho-of-derived-do} The cohomology of $\CD_P$ satisfy:
\[
H^i_\partial(\CD_P)=\left\{
\begin{array}{ccc} 0&\text{ if }&i<0,\\
\CD_A&\text{ if }&i=0.
\end{array}
\right.
\]
\end{lem}

Let $x_1,x_2,...$ be generators of degree 0 and $\xi_1,\xi_2,...$  generators of negative degree s.t. $P=\BC[x]\otimes\BC[\xi]$.
It is easy to see that $(\CD_P,[\partial,.])$ has a structure of a bicomplex s.t.  $\xi_i\mapsto (0,-1)$,
$\partial_{\xi_i}\mapsto (1,0)$.  The vertical differential creates a bunch of Koszul complexes, one for each monomial $\partial_{\xi_{i_1}}\partial_{\xi_{i_2}}\cdots$.
Therefore the spectral sequence degenerates at the term that has the form:
\begin{equation}
\label{2nd-term-spectr-seq}
0\rightarrow A\otimes_{\BC[x]}\CD_{\BC[x]}\rightarrow A\otimes_{\BC[x]}\CD_{\BC[x]}[\partial_{\xi_1},\partial_{\xi_2},\ldots]^{(1)}\rightarrow 
A\otimes_{\BC[x]}\CD_{\BC[x]}[\partial_{\xi_1},\partial_{\xi_2},\ldots]^{(2)}\rightarrow\cdots
\end{equation}
with differential $[\partial,.]$; here $A\otimes_{\BC[x]}\CD_{\BC[x]}[\partial_{\xi_1},\partial_{\xi_2},\ldots]^{(n)}$ means the space of polynomials in $\partial_{\xi}$'s of (cohomological) degree $n$
with coefficients in $A\otimes_{\BC[x]}\CD_{\BC[x]}$.
We leave it for the reader to convince himself that the 0th cohomology of this complex is $\CD_A$ essentially by definition.  $\qed$

It is natural therefore to think of $\CD_P$ as a resolution and thus a reasonable replacement of a  potentially pathological $\CD_A$. Of course, $\CD_P$ is determined
not so much by $A$ as by $P$, but one has:

\begin{lem} 
\label{homo-equi-derived-diff-oper}
The assignment $A\mapsto \CD_P$ defines an endofunctor on the homotopy category of commutative algebras and isomorphisms.
\end{lem}

{\em Proof.}  
V.Hinich (\cite{Hin} 8.5.3) proves that $A\mapsto \text{Der}(P)$ is a functor from the homotopy category of commutative algebras and isomorphisms to one of Lie algebras and isomorphisms; hence
for any other resolution $Q\rightarrow A$ there is a natural isomorphism (in the homotopy category) $\text{Der}(P)\rightarrow\text{Der}(Q)$. The algebras of differential operators carry a canonical filtration
s.t. the graded object is a symmetric algebra generated by derivations.  It suffices to prove that the arising morphism
\[
S^\bullet_P\text{Der}(P)\rightarrow S^\bullet_Q\text{Der}(Q)
\]
is an isomorphism.
Since $\text{Der}(P)$ and $\text{Der}(Q)$ are semi-free, their respective tensor powers (over $P$ and $Q$ resp.) are isomorphic, \cite{Hin} 3.3.2. In characteristic 0 this implies the desired isomorphism of symmetric powers.
$\qed$

\subsection{A chiral version.}  
\label{A chiral version.}
Unlike $\CD_A$, which is defined for any commutative ring $A$ (Grothendieck), its vertex algebra counterpart, $\CD^{ch}_{A,\vec{x}}$, introduced in 
sect.~\ref{local-defn-of-cdo} in the presence of a coordinate system,  has no obvious definition unless $\text{Spec}(A)$ is smooth. An obvious suggestion then is to
``chiralize'' the cohomological approach of sect.~\ref{Derived algebras of differential operators.}. Given $P\rightarrow A$ as above, we consider $\CD^{ch}_{P}$, which is allowed
as $P$ is a polynomial ring, and we drop the superscript $\vec{x}$ from the notation, as the coordinate system will be fixed once and for all.

The differential $\partial=\sum_if_i\partial_{x_i}$ has an obvious analogue: $\partial^{ch}=\sum_i(f_{i(-1)}\partial_{x_i})_{(0)}$. Note that although the product
$_{(-1)}$ is not commutative, the order does not matter. Indeed,
\[
f_{(i(-1)}\partial_{x_i}=\partial_{x_i(-1)}f_i-T(\partial_{x_i}(f_i)).
\]
Since $\text{deg}f_i=1$,  $\text{deg}\partial_{x_i}(f_i)$ is also 1. But by assumption $P$ is nonpositively graded, hence $\text{deg}\partial_{x_i}(f_i)\leq 0$. Therefore the
correction term $T(\partial_{x_i}(f_i))$ must be 0.

The element $\partial^{ch}\in\text{Der}\CD^{ch}_P$, see sect.~\ref{list-vert-alg-notat}. It is also a differential, as the following standard computation shows:
\[
2(\partial^{ch})^2=[\partial^{ch},\partial^{ch}]=[\partial,\partial]_{(0)}+\sum_{ij}\pm T(\partial_{x_j}(f_i))\partial_{x_i}(f_j).
\]
Here $[\partial,\partial]$ vanishes by assumption and the correction terms $T(\partial_{x_j}(f_i))\partial_{x_i}(f_j)$ vanish because, again, its degree must be 2 on the one
hand and $\leq 0$ on the other. Hence $(\partial^{ch})^2=0$.

We would like to think of the graded vertex algebra $(\CD^{ch}_P,\partial^{ch})$ as attached to $A$. Unfortunately, a satisfactory analogue of 
Lemma~\ref{homo-equi-derived-diff-oper} is unknown to us; what we do know will be explained later. Now we will change our point of view and talk about representation theory.

\subsection{Two categories of modules.}
\label{Two categories of modules.} 

\subsubsection{ }
\label{intro-to-vert-modules}
Let $V$ be a vertex algebra. A $V$-module is a vector space $\CM$ that carries a family of maps
\[
_{(j)}: V\otimes\CM\rightarrow\CM,\; v\otimes m\mapsto v_{(j)}m;\;j\in\BZ
\]
s.t. a number of axioms are satisfied; see e.g. \cite{FBZ} for details.

In the case where $V$ is conformally graded, $V=\oplus_{n\geq 0}V[n]$, a $V$-module $\CM$ is called conformally graded if there is a direct sum decomposition
$\CM=\oplus_{m\in \BZ}\CM[m]$ s.t. $V[n]_{(j)}\CM[m]\subset\CM[n+m-j-1]$, cf. sect.~\ref{list-vert-alg-notat}. 

Finally, if $V$ is a differential graded vertex algebra with differential $\partial^{ch}$, then a differential graded $V$-module is a graded $V$-module $\CM$ with degree 1 differential
$\partial_{\CM}:\CM\rightarrow\CM$ s.t. $\partial_{\CM}(v_{(j)}m)=\partial^{ch}(v)_{(j)}m+(-1)^{\text{deg}v}v_{(j)}\partial_{\CM}(m)$ for all $j$. The grading thus built into
the definition is sometimes referred to as cohomological grading. If, in addition, $V$ and $\CM$ are conformally graded, then $\partial_\CM$ is required to have conformal
weight 0.

\subsubsection{ } 
\label{usual dp--into play}
Representation theory is where $\CD_P$ and $\CD^{ch}_P$ actually meet.  One easily verifies the relations (among elements of $\text{End}\CD^{ch}_P$):
\begin{equation}
\label{comm-rel-x-dx}
[\partial_{x_i(m)},x_{j(n-1)}]=\delta_{i,j}\delta_{m,-n}Id.
\end{equation}
This implies that for any $\CD^{ch}_P$-module $\CM$, the elements $\partial_{x_i(m)},x_{j(n-1)}$, regarded as operators acting on $\CM$, satisfy the same relations.
Since $[\partial_{x_i(0)},x_{j(-1)}]=\delta_{i,j}Id$, which is a familiar bracket of  coordinate vector fields and functions, every $\CD^{ch}_P$-module becomes a
$\CD_P$-module.  Since   $x_i\in\CD^{ch}_P[0]$,  $\partial_{x_i}\in\CD^{ch}_P[1]$, cf. the last sentence of sect.~\ref{local-defn-of-cdo}, 
the corresponding operators $x_{j(-1)}$ and $\partial_{x_i(0)}$ have conformal weight 0; therefore each homogeneous component
$\CM[m]$ is a $\CD_P$-submodule of $\CM$.

Consider the space of singular vectors
\[
\text{Sing}\CM=\{m\in\CM\text{ s.t. }\partial_{x_i(n)}m=x_{j(n)}m=0\;\forall n>0.\}.
\]
Since the linear span of $\{x_{i(n)},\partial_{x_j(n)},\;n>0\}$ is preserved under the bracket with $\{x_{i(-1)},\partial_{x_j(0)},\}$, $\text{Sing}\CM \subset\CM$ is a 
$\CD_P$-submodule. Furthermore, if in addition $\CM$ is a dg-module, the restriction of $\partial^{ch}_\CM$ to $\text{Sing}\CM$ endows the latter with a structure
of a dg-module over $(\CD_P,\partial)$. To summarize, the rule $\CM\mapsto\text{Sing}\CM$ defines a functor from the category of dg-modules over $\CD^{ch}_P$
to the category of  dg-modules over $\CD_P$.

There is no reason to believe that $\text{Sing}\CM$ is in general nonzero. To ensure that it is, let us make the following definition: call $\CM$ {\em bounded} if for each $m\in\CM$
there is $N>0$ s.t. $x_{i(n)}^Nm=\partial_{x_i(n)}^Nm=0$ for all $i$ and $n>0$. 

Note that for each  $m$, $x_{i(n)} m=\partial_{x_i(n)} m=0$ for all sufficiently large $n$ -- this is
part of the definition of a vertex algebra module..
Therefore boundedness simply  means that the elements $x_{i(n)},\partial_{x_i(n)}$, $n>0$, act locally nilpotently.

Of course, if $\CM$ is graded and $\CM[n]=\{0\}$ if $n\gg 0$, then $\CM$ is bounded; in
particular, $\CD^{ch}_P$ is bounded if considered as a module over itself.  However, the definition does not require any grading.

\begin{lem}
\label{equi-dg-dch-mod-dg-d-mod}
The assignment $\CM\mapsto\text{Sing}\CM$ defines an equivalence of  the category of bounded dg-$\CD^{ch}_P$-modules and the category of
dg-$\CD_P$-modules.
\end{lem}

{\em Proof.} To begin with, let us forget about the dg-structure; in this case we have: {\em $\CM\mapsto\text{Sing}\CM$ defines an equivalence of the category of
bounded $\CD^{ch}_P$-modules and the category of $\CD_P$-modules.}

This result was proved by D.Chebotarov \cite{Ch}, but since he worked in a much more general and complicated situation, we will give an independent proof.
Denote by $Mod(\CD_P)$ and $Mod(\CD^{ch}_P)$ the categories in question and let $\Phi:Mod(\CD^{ch}_P)\rightarrow Mod(\CD_P)$ be the functor we have defined.
We need a functor in the opposite direction: $\Psi:Mod(\CD_P)\rightarrow Mod(\CD^{ch}_P)$, and we will show that  $\Psi$ can be identified with the 
 induction functor for a certain infinite dimensional Lie algebra.

Let $\hat{\fa}$ denote the Lie algebra  with basis $\{a^i_n,b^i_n, C;\; n\in\BZ,1\leq i\leq \text{dim}P\}$ and relations
\[
[a^i_n,b^j_m]=\delta_{i,j}\delta_{n,-m}C,\; [C,a^i_n]=[C,a^j_m]=0.
\]
A glance at \ref{comm-rel-x-dx} will show that the assignment $a^i_n\mapsto\partial_{x_i(n)}$, $b^i_n\mapsto x_{i(n-1)}$, $C\mapsto Id$
makes any $\CD^{ch}_P$-module into an $\hat{\fa}$-module. Furthermore, one defines what a bounded $\hat{\fa}$-module is in an obvious manner (and requires
that $C$ act as $Id$) and thus
obtains a functor on the categories of bounded modules: $Mod(\CD^{ch}_P)\rightarrow Mod(\hat{\fa})$. It is a standard result that this functor is an equivalence
of categories; see e.g. \cite{FBZ}, 5.1.8. From now on we will not distinguish between these 2 categories.

Notice that $\hat{\fa}$ is graded: $\hat{\fa}=\oplus_n\hat{\fa}[n]$, where $\hat{\fa}[n]$, $n\not=0$, is spanned by $a^i_{-n}$ and $b^j_{-n}$, and $\hat{\fa}[0]$ is spanned by
$a^i_0$, $b^j_0$, and $C$. It is obvious that $U(\hat{\fa}[0])/\langle C-1\rangle =\CD_P$.   Denote by $\hat{\fa}_>=\oplus_{n>0}\hat{\fa}[n]$,
 $\hat{\fa}_\geq=\oplus_{n\geq 0}\hat{\fa}[n]$. There arises a Lie algebra homomorphism $U(\hat{\fa}_\geq)\rightarrow\CD_P$ s.t. $\hat{\fg}_>\mapsto 0$.
 Therefore any $\CD_P$-module, by pull-back, becomes a $\hat{\fa}_\geq$-module. Now define
 \[
 \Psi\stackrel{\text{def}}{=}\text{Ind}_{\hat{\fa}_\geq}^{\hat{\fa}}: Mod(\CD_P)\rightarrow Mod(\CD^{ch}_P).
 \]
 By definition, $\Phi\circ\Psi: Mod(\CD_P)\rightarrow Mod(\CD_P)$ is the identity. To complete the proof we need show that $\Psi\circ\Phi$ is isomorphic to the identity
 functor. At least there is a functor morphism $\Psi\circ\Phi\rightarrow Id$ as for each $\CM\in Mod(\CD^{ch}_P)$ there is a canonical morphism
 \[
 \text{Ind}_{\hat{\fa}_\geq}^{\hat{\fa}}\text{Sing}\CM\rightarrow\CM
 \]
 which is determined (thanks to the universal property of the induction) by the tautological inclusion $\text{Sing}\CM\hookrightarrow\CM$. We need to show that this map is
 injective and surjective. Both these assertions follow easily from the properties  of the grading operator
 \[
 L_0\stackrel{\text{def}}{=}\sum_{i=1}^{\text{dim}}\sum_{n=1}^{+\infty}-n(b^i_n a^i_{-n}-a^i_nb^i_{-n}).
 \]
Note that the eigenvalues of $[L_0,.]$ define the above grading of $\hat{\fa}$: $[L_0,a^i_n]=-na^i_n$, $[L_0,b^i_n]=-nb^i_n$. It follows that $ \text{Ind}_{\hat{\fa}_\geq}^{\hat{\fa}}\text{Sing}\CM$ is also graded by the eigenvalues of $L_0$:  
\[
 \text{Ind}_{\hat{\fa}_\geq}^{\hat{\fa}}\text{Sing}\CM=\oplus_{n\geq 0} \text{Ind}_{\hat{\fa}_\geq}^{\hat{\fa}}\text{Sing}\CM[n].
 \]
 Furthermore $ \text{Ind}_{\hat{\fa}_\geq}^{\hat{\fa}}\text{Sing}\CM[0]=\text{Sing}\CM$, and the rest of the eigenvalues are strictly positive.
 
 {\em Injectivity.}  The kernel of the morphism $\text{Ind}_{\hat{\fa}_\geq}^{\hat{\fa}}\text{Sing}\CM\rightarrow\CM$ is a $\hat{\fa}$-submodule. Therefore the kernel must
 contain a singular vector, i.e., a vector $v$ that satisfies $a^i_nv=b^i_nv=0$ for all $n>0$. This forces $L_0v=0$, hence $v\in\text{Sing}\CM$, therefore $v=0$.
 
 {\em Surjectivity.} Denote by $\CM'$ the $\hat{\fa}$-submodule of $\CM$ generated by $\text{Sing}\CM$. By definition, the image of  the map
$\text{Ind}_{\hat{\fa}_\geq}^{\hat{\fa}}\text{Sing}\CM\rightarrow\CM $ is precisely $\CM'$. Let $\overline{\CM}$ be $\CM/\CM'$. We need to show that $\overline{\CM}=\{0\}$.

It is clear from the definition of a bounded module that $\overline{\CM}$ is bounded and, therefore, contains a singular vector $\bar{v}$. As above, $L_0\bar{v}=0$. 
Let $v\in\CM$ be a preimage of $\bar{v}$.  It is clear that $v=v_1+v_2$, where $v_1\in\CM'$ and $v_2$ is a generalized eigenvector of eigenvalue 0, i.e.,
$L_0^Nv_2=0$ for some $N$.  Hence $a^i_nv_2$, $b^i_nv_2$ are generalized eigenvectors of eigenvalue $-n$. Since $\bar{v}$ is singular
 these vectors belong to $\CM'$ for $n>0$ where, as
we have seen, all eigenvalues of $L_0$ are nonnegative. This means that $a^i_nv_2=b^i_nv_2=0$ if $n>0$. Therefore $v_2\in\text{Sing}\CM$ and $\bar{v}=0$, as desired.

Finally, we have to decide what $\Psi$ does to the differential graded structure. As we have seen, any $\CD^{ch}_P$-module $\CM$ is induced from its space of singular
vectors $\text{Sing}\CM$.  This implies at once that a dg-$\CD_P$-module structure on the latter has at most one extension to a dg-$\CD^{ch}_P$-module structure on $\CM$.
It suffice to prove the existence. 

If, as in our situation, the differential on the algebra is defined by an algebra element ($[\partial,.]$, $\partial\in\CD_P$ in the case of $\CD_P$, and $\partial^{ch}_{(0)}$,
$\partial^{ch}\in\CD^{ch}_P$ in the case of $\CD^{ch}_P$),  the differential on a dg-module is uniquely written in the form: the corresponding element of the algebra plus
a square 0
 endomorphism of the module of cohomological degree 1. For example, the differential on $\text{Sing}\CM$ equals $\partial+f$, where $f\in\text{End}_{\CD_P}(\text{Sing}\CM)$. The corresponding differential on $\CM$ is then $\partial^{ch}_{(0)}+\Psi(f)$.
 
 This construction extends $\Psi$ to a functor between categories of dg-modules, $\Phi$ was previously defined on dg-modules, and it is quite clear that the established above
 isomorphisms of functors  $\Phi\circ\Psi = Id$ and $\Psi\circ\Phi\iso Id$ respect the dg-module structure. $\qed$
 
 \begin{cor}
 \label{mor-equi-of-chir-res}
 For any two finite polynomial resolutions $P\rightarrow A\leftarrow Q$ the categories of dg modules over $\CD^{ch}_P$ and $\CD^{ch}_Q$ are
 equivalent.
 \end{cor}
 
 \begin{cor}
 \label{any-bounded-graded}
 Any bounded $\CD^{ch}_P$-module $\CM$ can be conformally  graded so that $\CM[0]=\text{Sing}\CM$ and $\CM[n]=0$ if $n<0$.
 \end{cor}

\subsubsection{Example.}  
\label{Example.}
Let $A=\BC[x]/(x^n)$.  As $P$ take the Koszul resolution w.r.t. the regular sequence consisting of one element $x^n$; this is what in
sect.~\ref{derived-version-of-a-ring} was denoted by $K(\BC[x],x^n)$. This is a polynomial ring on 2 variables, $x$, which is even, $\xi$ odd (of degree -1),
best thought of as the struture ring of  the superspace $\BC^{1|1}$.
The corresponding dg algebra of differential operators is $(\CD_{\BC^{1|1}},[x^n\partial_\xi,.])$.

\begin{lem}
\label{do-fat-point-quis-gln}
The dg algebra $(\CD_{\BC^{1|1}},[x^n\partial_\xi,.])$ is quasiisomorphic to $gl_n(\BC)$.
\end{lem}

{\em Proof.} In this case the complex (\ref{2nd-term-spectr-seq}) becomes
\[
0\rightarrow\CD_\BC/x^n\CD_\BC\stackrel{[x^n\partial_\xi,.]}{\longrightarrow}\partial_\xi\CD_\BC/x^n\CD_\BC\rightarrow 0.
\]
We will show that the differential is surjective, and its kernel  has dimension $n^2$.  Since by virtue of Lemma~\ref{coho-of-derived-do} 
$H^0_\partial(\CD_{\BC^{1|1}})=\CD_{\BC[x]/(x^n)}\hookrightarrow gl(\BC[x]/(x^n))$, which is $n^2$-dimensional, the lemma will follow.

Proving these 2 assertions requires doing  a bit of  linear algebra. The space $\CD_\BC$ is graded by $x\mapsto 1$, $\partial_x\mapsto -1$. If we denote by $\CD_{\BC,j}$
the homogeneous component of degree $j$, then $\CD_{\BC,j}=\{0\}$ if $j\geq n$, $\text{dim}\CD_{\BC,j}=n-j$ if $0\leq j\leq n-1$; thus 
$\text{dim}\CD_{\BC,0}=n$ at which point the dimension stabilizes s.t. $\text{dim}\CD_{\BC,j}=n$ if $j<0$.

The differential is apparently of degree $n$. Therefore the components $\CD_{\BC,j}$ with $j=n-1,n-2,...,0$ belong to the kernel in their entirety by dimensional argument.
This contributes $1+2+3+\cdots n$ to the dimension of the kernel.

For the next $n-1$ components, the differential gives us maps
\[
\CD_{\BC,j}\rightarrow\CD_{\BC,n+j},\; j=-1,-2,...,-n+1.
\]
We will show that all these maps are surjections. Since the $j$-th such map has $n$-dimensional domain and $(-j)$-dimensional codomain, its kernel is $(n+j)$-dimensional.
This contributes $(n-1)+(n-2)+\cdots +1$ to the dimension.

The remaining maps
\[
\CD_{\BC,j}\rightarrow\CD_{\BC,n+j},\; j=-n,-n-1,....
\]
are also surjective, but have domain and codomain of the same dimension $n$, hence are isomorphisms.

Overall the dimension of the kernel is $1+2+\cdots+(n-1)+n+(n-1)+\cdots+2+1=n^2$ as desired.

To see that the maps $\CD_{\BC,j}\rightarrow\CD_{\BC,n+j},\; j<-1$ are surjective, one needs to write down their matrices. One has
\[
[x^n,\partial_x^j]=\sum_{i=1}^{\text{min}\{n,j\}}\alpha_ix^{n-i}\partial_x^{j-i},
\]
where none of the $\alpha$'s vanishes. From this the untiring reader will deduce that when written in the basis of the form 
$\{\partial_x^k,x\partial_x^{k+1},\ldots x^{n-1}\partial_x^{k+n-1}\}$ the images of the last $k$ (or all $n$ if $k>n$) vectors form a triangular matrix with nonzero
diagonal entries. This completes the proof. $\qed$

\begin{cor}
\label{coroll-for-cdo-fat-point}
The derived category of dg-modules over the dg vertex algebra $(\CD^{ch}_{\BC^{1|1}})$ is equivalent to the derived category of finite dimensional vector
spaces.
\end{cor}

{\em Proof.} Quasiisomorphic dg algebras have equivalent derived dg module categories; this is \cite{Hin} Theorem~3.3.1. Now Corollary~\ref{coroll-for-cdo-fat-point} follows
at once from Lemma~\ref{equi-dg-dch-mod-dg-d-mod} and the classic result that  the module category of a matrix algebra is semi-simple on one generator (one can consider this particular case of the Morita 
equivalence as an "odd Stone - von Neumann theorem"). $\qed$

This corollary has a pleasing interpretation.
\subsubsection{ }
\label{interpr-feig-sem}
The Landau-Ginzburg model  starts with a function (potential) defined on $\BC^n$ and is ultimately related to the $N=$ superconformal
algebra \cite{Witt}.   If the potential $\Phi\in\BC[x_1,...,x_n]$, one can consider the above defined $\CD^{ch}_P$ 
with $P=K(\BC[x_1,...,x_n],\{\Phi'_{x_1},...,\Phi'_{x_n}\})$, the Koszul resolution of $A=\BC[x_1,...,x_n]/(\Phi'_{x_1},...,\Phi'_{x_n})$.
The cohomology of this algebra was computed by Feigin and Semikhatov \cite{FS} in the case where $n=1$, $\Phi=x^{N+1}$, and their answer
is a thing of beauty: the cohomology is a direct sum of $N$ vector spaces, one vector space for each of the classical Koszul classes represented by
$1,x,...,x^{N-1}$; the space attached to 1 is the unitary vacuum representation of the $N=2$ superconformal algebra, and the rest of the spaces are
irreducible unitary representations of this algebra.

Of course this result immediately carries over to the case of a ``diagonal'' potential, such as $\Phi=\sum_ix_i^n$, by using the K\"unneth formula.
For this reason, we will denote by $\CF\CS_\Phi$ the dg vertex algebra $\CD^{ch}_P$ 
with $P=K(\BC[x_1,...,x_n],\{\Phi'_{x_1},...,\Phi'_{x_n}\})$ .

A vertex algebra is called {\em conformal} if it contains a conformal vector $L$, i.e., one that ``generates '' the Virasoro algebra with some central charge,
see \cite{FBZ} 2.5.8. The vertex algebra $\CD_P^{ch}$ is conformal with conformal  vector $L=Tx_{i(-1)}\partial_{x_i}$, cf.\cite{MSV}. 

We will call a dg vertex algebra {\em dg conformal} if it is conformal and the conformal vector is annihilated by the differential. It is easy to verify that 
$\partial^{ch}(Tx_{i(-1)}\partial_{x_i})=0$, and so $(\CD^{ch}_P,\partial^{ch})$ is dg conformal.

A vertex algebra is called {\em rational}, see \cite{FBZ} 5.5.1 and references therein, if it is conformal and  its module category
is semi-simple. It is then natural to call a dg vertex algebra {\em derived rational} if the derived category of its dg module category is equivalent to
the derived category of a semi-simple category.

Corollary~\ref{coroll-for-cdo-fat-point} gives:
\begin{thm}
\label{f-sem-vert-alg-ratnl}
If $\Phi$ is diagonal, then the Feigin-Semikhatov vertex algebra $\CF\CS_\Phi$ is derived rational.
\end{thm}

\subsubsection{}
As was explained to us by D.Gaitsgory, it follows from \cite{GR} that
the derived category of dg $\CD_P$-modules is equivalent to the derived category of $\CD$-modules on $\text{Spec}(A)$.
In conjunction with Lemma~\ref{equi-dg-dch-mod-dg-d-mod} this implies that the derived category of dg-$\CD^{ch}_P$-modules is equivalent to the latter category. This
is a generalization of Corollary~\ref{coroll-for-cdo-fat-point}, because the category of $\CD$-modules on a fat point, which is by definition
 the category of $\CD_\BC$-modules
supported on the corresponding closed point, is isomorphic to the category of finite dimensional vector spaces. More generally, the Gaitsgory's
 result implies that
$\CD^{ch}_P$ is derived rational provided  $P$  is a resolution of the structure ring of a finite collection of fat points; this, of course, includes the case
of $\CF\CS_\Phi$ if $\Phi$ is a polynomial with isolated singularities, not necessarily diagonal. 

\subsection{The graded ring case}
\label{the-graded-ring-case}

We will now explain what we know about the naturality of the assignment $A\mapsto\CD^{ch}_P$.

\subsubsection{ } We wish to add a 3rd type of grading to the two, conformal and cohomological, that have prominently figured
so far. 
By an {\em inner} grading on an algebra $A$  we will mean a vector space decomposition
 $A=A_0\oplus A_{1}\oplus A_{2}\oplus\cdots$ s.t. $A_i A_j\subset A_{i+j}$ and
$A_0=\BC1$. $M$ is a graded $A$-module if similarly $M=M_0\oplus M_1\oplus
\cdots$ with $A_i M_j\subset M_{i+j}$.

In this section the phrase ``$A$ is a graded algebra'' will mean an algebra $A$ carrying an inner
grading.

A grading of a polynomial ring $\BC[x_1,...,x_n]$ is the same as an assignment of a positive integer $\text{deg}\ x_i$ to each $x_i$.

Call a graded algebra $A$ a complete intersection  if it is a quotient of a graded polynomial ring 
$\BC[x_1,...,x_n]$ by a homogeneous
ideal $J$ that is generated by a homogeneous regular sequence. Having fixed one such sequence $\vec{f}=\{f_1,...,f_m\}$ we obtain a Koszul complex $K(\BC[\vec{x},\vec{\xi}],\vec{f})$
and a dga quasiisomorphism $K(\BC[\vec{x},\vec{\xi}],\vec{f})\rightarrow A$. As  algebra, $K(\BC[\vec{x},\vec{\xi}],\vec{f})=\BC[x_1,...,x_n,\xi_1,...\xi_m]$, $\xi_i$ having cohomological degree -1
and inner degree $\text{deg}f_i$, $i=1,2,...$; the differential is $\partial=\sum f_i\partial_{\xi_i}$.

This gives us a differential graded vertex algebra $D^{ch}(K(\BC[\vec{x},\vec{\xi}],\vec{f}))$, sect.~\ref{A chiral version.}.

Call two differential graded vertex algebras $V$ and $W$ {\em quasiisomorphic} if there is a sequence of differential graded vertex algebra quasiisomorphisms
\[
 V\leftarrow V_0\rightarrow V_1\leftarrow\cdots\rightarrow W.
\]

\begin{thm}
 \label{result-graded-case}
For any two quasiisomorphisms $K(\BC[x_1,...,x_n],\vec{f})\rightarrow A\leftarrow K(\BC[y_1,...,y_l],\vec{g})$
the chiral Koszul complexes $D^{ch}(K(\BC[\vec{x},\vec{\xi}],\vec{f}))$ and $D^{ch}(K(\BC[\vec{y},\vec{\eta}],\vec{g})$ are
quasiisomorphic.
\end{thm}
The proof will occupy a few subsections to follow.

\subsubsection{  }
\label{minimal-gener-set}
Let $M$ be a {\em finitely} graded module over a graded algebra $A$.  Let $A_>=\oplus_{n\geq 1}A_n$.
A {\em minimal generating space} of $M$ is the image of a homogeneous splitting
of the projection
\[
M\longrightarrow M/A_>M.
\]
A {\em minimal generating set} of $M$ is a choice of a homogeneous basis of a minimal generating space.

Any 2 minimal generating sets, $E$ and $F$, are related by a transformation matrix $a_{E\rightarrow F}$ with coefficients in $A$.
$M$ being not necessarily free, $a_{E\rightarrow F}$ is not unique, but it can always be chosen so as to be a block upper triangular, diagonal
blocks being constant matrices. This implies that such  $a_{E\rightarrow F}$ is automatically invertible.

\subsubsection{ } \label{lifting-coord-changes}
 Fix a Koszul complex $K(\BC[\vec{x},\vec{\xi}],\vec{f})$.  As a free $\BC[\vec{x}]$-module, it carries, according to \ref{minimal-gener-set}, an action of the group
 of coordinate transformations. Since we would like to have the cohomological degree preserved, any such coordinate transformation will be a composition of the transformations
  of the
 following 2 types:
 \[
 x_i\mapsto f_i(x)\text{ or }\xi_i\mapsto a_{ij}(x)\xi_j.
 \]
 We claim that these coordinate changes can be lifted (non-canonically) to automorphisms of the corresponding chiral Koszul complex $D^{ch}(K(\BC[\vec{x},\vec{\xi}],\vec{f}))$.
 Indeed,  the coordinate change of type 1 classically gives rise to
 \[
 \partial_{x_i}\mapsto\sum_{s}\frac{\partial g_s(x)}{\partial x_i}\partial_{x_s},\;\xi\mapsto\xi_i,\partial_{\xi_i}\mapsto\partial_{\xi_i},\; g=f^{-1}.
 \]
Define a vertex algebra lift  to be simply
 \[
 \partial_{x_i}\mapsto\sum_{s}(\partial_{x_s})_{(-1)}(\frac{\partial g_s(x)}{\partial x_i}),\;\xi\mapsto\xi_i,\partial_{\xi_i}\mapsto\partial_{\xi_i}.
 \]
It defines an automorphism, because by the Wick theorem (cf. sect.~\ref{list-vert-alg-notat})
 \[
(\partial_{x_s})_{(-1)}(\frac{\partial g_s(x)}{\partial x_i})(z)\cdot (\partial_{x_s})_{(-1)}(\frac{\partial g_s(x)}{\partial x_j})(w)=
\frac{-1}{(z-w)^2}\frac{\partial^2 g_t(x)}{\partial x_s\partial x_j}(z)\frac{\partial^2 g_s(x)}{\partial x_t\partial x_i}(w)+\cdots,
\]
which is 0, as it should,  by virtue of the following obvious dimensional argument: the degree of the expression
$\frac{\partial^2 g_t(x)}{\partial x_s\partial x_j}(z)\frac{\partial^2 g_s(x)}{\partial x_s\partial x_i}(w)$ is $-deg x_i-deg x_j$, hence strictly negative, while the expression itself is a product
of fields coming from elements of $\BC[\vec{x}]$, $\frac{\partial^2 g_t(x)}{\partial x_s\partial x_j}$ and
$\frac{\partial^2 g_s(x)}{\partial x_s\partial x_i}$, hence has nonnegative degree.

Now consider the type 2 transformation $\xi_i\mapsto a_{ij}(x)\xi_j$.  An obvious vertex algebra lift is defined by
\[
\partial_{x_i}\mapsto \tilde{\partial}_{x_i}=\partial_{x_i}+(\partial_{x_i}b_{sj})a_{jt}\xi_t\partial_{\xi_s},
\; \partial_{\xi_i}\mapsto b_{ji}\partial_{\xi_j},
\]
where $(b_{ij})$ is the inverse of $(a_{ij})$ -- which exists by \ref{minimal-gener-set}. Another 
application of the Wick theorem gives
\[
\tilde{\partial}_{x_i}(z)\tilde{\partial}_{x_j}(w)=Tr(\partial_{x_i}b)a(z)(\partial_{x_j}b)a(w)+\cdots.
\]
We know from \ref{minimal-gener-set} that in an appropriate basis $(a_{ij})$ is block upper-triangular with constant diagonal blocks. Therefore
$Tr(\partial_{x_i}b)a(z)(\partial_{x_j}b)a(w)=0$, as desired.

\subsubsection{ }
\label{minimal-koszul-res}
Consider $A_>\subset A$ as a graded $A$-module and
pick a minimal generating set $\bar{X}=\{\bar{x}_1,...,\bar{x}_n\}\subset\CA$. This gives an algebra surjection $\BC[x_1,...,x_n]\rightarrow\CA$ with kernel $J$.
Pick a minimal generating set $\{f_1,...,f_m\}\subset J$, where $J$ is regarded as a $\BC[x_1,...,x_n]$-module. This defines  a Koszul complex $K(\BC[\vec{x},\vec{\xi}],\vec{f})$ and its chiral version
$D^{ch}(K(\BC[\vec{x},\vec{\xi}],\vec{f}))$. We will refer to these as a {\em minimal} (chiral) Koszul complex.

We assert that all minimal, chiral or usual, Koszul complexes are isomorphic.  Indeed,  for any other minimal generating set $\bar{Y}=\{\bar{y}_1,\bar{y}_2,...\}$, there
is a transformation matrix $a_{X\rightarrow Y}$, which can be lifted to $\BC[\vec{y}]$. Any such lift  defines an algebra isomorphism $\BC[\vec{x}]\rightarrow\BC[\vec{y}]$
that respects both projections $\BC[\vec{x}]\rightarrow A\leftarrow\BC[\vec{y}]$.

Next, consider the Koszul complex $K(\BC[\vec{y},\vec{\eta}],\vec{g})$ corresponding to $\bar{Y}$. Under the constructed isomorphism $\BC[\vec{x}]\rightarrow\BC[\vec{y}]$,
each $f_i$ is mapped to $a_{ij}g_j$ for some  transformation matrix $(a_{ij})$. Composing with $\xi_i\mapsto  a_{ij}\eta_j$ gives us an isomorphism of complexes
$K(\BC[\vec{x},\vec{\xi}],\vec{f})\rightarrow K(\BC[\vec{y},\vec{\eta}],\vec{g})$.

It follows from  \ref{lifting-coord-changes} that both these  transformations  lift to the chiral Koszul complex and their composition defines the desired isomorphism
\[
D^{ch}(K(\BC[\vec{x},\vec{\xi}],\vec{f}))\rightarrow D^{ch}(K(\BC[\vec{y},\vec{\eta}],\vec{g})).
\]
\subsubsection{ }
\label{from-any-to-minimal}
We will now show that for any Koszul resolution $K(\BC[\vec{x},\vec{\xi}],\vec{f})\rightarrow A$ there is a minimal Koszul resolution
$K_{min}\rightarrow A$ and a quasiisomorphism $D^{ch}(K(\BC[\vec{x},\vec{\xi}],\vec{f}))\rightarrow D^{ch}(K_{min})$.

The elements of the regular sequence $\vec{f}=\{f_1,f_2,...,f_m\}$ break down into 2 groups: those that  do not belong to $\BC[\vec{x}]_{>}^2$ and those that do.
The former group is not necessarily linearly independent modulo  $\BC[\vec{x}]_{>}^2$, but it is clear that replacing $\{f_1,f_2,...,f_m\}$ with their appropriate
linear combinations   we ensure that it is.

 If so, upon an appropriate coordinate transformation this  part of $\vec{f}$ becomes part of a coordinate system.  We can therefore assume that
 $\vec{f}=\{x_1,x_2,...,x_k,f_{k+1},...,f_n\}$.  Let $J=(x_1,...,x_k,\xi_1,...,\xi_k)\subset \BC[\vec{x},\vec{\xi}]$ and define $K_{min}=\BC[\vec{x},\vec{\xi}]/J$.
 $J$ is a differential graded ideal, hence  $\BC[\vec{x},\vec{\xi}]/J$ is a differential graded algebra; in fact it is a Koszul complex with differential
 $\partial=\sum_{i>k}\bar{f}_{i}\partial_{\bar{\xi}_i}$. Furthermore, $J$ is acyclic ($\xi_i\mapsto x_i$, $x_i\mapsto 0$) and so is $K_{min}$. This implies that $K_{min}$ is a minimal Koszul
 resolution of  $A$, and we have a quasiisomorphism
 \[
K(\BC[\vec{x},\vec{\xi}],\vec{f})\rightarrow K_{min}.
\]
All of this ``chiralizes.'' We start with $D^{ch}(K(\BC[\vec{x},\vec{\xi}],\vec{f}))$.  Then we lift to $D^{ch}(K(\BC[\vec{x},\vec{\xi}],\vec{f}))$, one-by-one, the 2 coordinate transformations
that include part of $\vec{f}$ into a coordinate system; this is based on \ref{lifting-coord-changes}. Finally, we quotient out a vertex ideal ideal, $J_{vert}$,  generated by
$x_1,...,x_k;\partial_{x_1},...,\partial_{x_k};\xi_1,...,\xi_k;\partial_{\xi_1},...,\partial_{\xi_k}$. The result is an exact sequence
\[
0\rightarrow J_{vert}\rightarrow D^{ch}(K(\BC[\vec{x},\vec{\xi}],\vec{f}))\rightarrow D^{ch}(K_{min})\rightarrow 0.
\]
$J_{vert}$ is acyclic (  $\xi_i\mapsto x_i$, $\partial_{x_i}\mapsto-\partial_{\xi_i}$.) Therefore, the projection
\[
D^{ch}(K(\BC[\vec{x},\vec{\xi}],\vec{f}))\rightarrow D^{ch}(K_{min})
\]
is a quasiisomorphism.

Since all minimal chiral Koszul resolutions are isomorphic, \ref{minimal-koszul-res}, for any 2 Koszul resolutions
\[
K(\BC[\vec{x},\vec{\xi}],\vec{f})\rightarrow A\leftarrow K(\BC[\vec{y},\vec{\eta}],\vec{g})
\]
we obtain a diagram  consisting of quasiisomorphisms
\[
D^{ch}(K(\BC[\vec{x},\vec{\xi}],\vec{f}))\rightarrow D^{ch}(K_{min})\leftarrow D^{ch}(K(\BC[\vec{y},\vec{\eta}],\vec{g})).
\]
This concludes the proof of Theorem\ref{result-graded-case}.

\subsubsection{ }
\label{rem-on-anal} A similar argument gives an analogous result for complete intersections
in  analytic setting without the grading assumption.  This is an evidence for there existing
a gerbe of dg vertex algebras over any locally complete intersection at least in analytic setting.

\bigskip

F.M.: Department of Mathematics, University of Southern California, Los Angeles, CA 90089, USA

V.Sch.: Institut de Math\'ematiques de Toulouse, Universit\'e Paul Sabatier, 31062 Toulouse, France

\end{document}